\documentclass[12pt,leqno]{amsart}
\usepackage{amssymb,verbatim,enumerate,ifthen,amsmath}
\usepackage[mathscr]{eucal}
\usepackage[utf8]{inputenc}
\usepackage[T1]{fontenc}
\def\N{\mathbb{N}}
\def\R{\mathbb{R}}

\def\F{\mathscr{F}}
\def\C{\mathscr{C}}

\def\H{\mathscr{H}}

\newtheorem{theorem}{Theorem}[section]
\newtheorem*{theorem*}{Theorem}
\def\Thm#1#2{\ifthenelse{\equal{#1}{*}}{\begin{theorem*}#2\end{theorem*}}
             {\begin{theorem}\label{T#1}#2\end{theorem}}}
\newtheorem{Atheorem}{Theorem}

\def\thm#1{Theorem~\ref{T#1}}
\newtheorem{proposition}[theorem]{Proposition}
\newtheorem*{proposition*}{Proposition}
\def\Prp#1#2{\ifthenelse{\equal{#1}{*}}{\begin{proposition*}#2\end{proposition*}}
{\begin{proposition}\label{P#1}#2\end{proposition}}}
\def\prp#1{Proposition~\ref{P#1}}

\newtheorem{corollary}[theorem]{Corollary}
\newtheorem*{corollary*}{Corollary}
\def\Cor#1#2{\ifthenelse{\equal{#1}{*}}{\begin{corollary*}#2\end{corollary*}}
             {\begin{corollary}\label{C#1}#2\end{corollary}}}

\newtheorem{lemma}[theorem]{Lemma}
\newtheorem*{lemma*}{Lemma}
\def\Lem#1#2{\ifthenelse{\equal{#1}{*}}{\begin{lemma*}#2\end{lemma*}}
             {\begin{lemma}\label{L#1}#2\end{lemma}}}

\theoremstyle{definition}
\newtheorem{remark}[theorem]{Remark}
\newtheorem*{remark*}{Remark}
\def\Rem#1#2{\ifthenelse{\equal{#1}{*}}{\begin{remark}\rm #2\end{remark}}
             {\begin{remark}\label{R#1}\rm #2\end{remark}}}

\newtheorem{example}[theorem]{Example}
\newtheorem*{example*}{Example}
\def\Exa#1#2{\ifthenelse{\equal{#1}{*}}{\begin{example*}\rm #2\end{example*}}
             {\begin{example}\label{Ex#1}\rm #2\end{example}}}

\def\eq#1{{\rm(\ref{E#1})}}
\def\Eq#1#2{\ifthenelse{\equal{#1}{*}}
  {\begin{equation*}\begin{aligned}#2\end{aligned}\end{equation*}}
  {\begin{equation}\begin{aligned}\label{E#1}#2\end{aligned}\end{equation}}}

\begin{document}
\begin{flushright}
\end{flushright}
\vspace{5mm}

\date{\today}

\title{Generalization of Subadditive, Monotone and Convex Functions}

\author[A. R. Goswami]{Angshuman R. Goswami}
\address{Department of Mathematics, University of Pannonia, 
Veszprém, Egyetem u. 10, H-8200.}
\email{goswami.angshuman.robin@mik.uni-pannon.hu}


\begin{abstract}
Let $I\subseteq\R_+$ be a non empty and non singleton interval where $\R_+$ denotes the set of all non negative numbers. A function $\Phi: I\to \R_+$ is said to be subadditive if for any $x,y$ and $x+y\in I$, it satisfies the following inequality 
$$\Phi(x+y)\leq \Phi(x)+\Phi(y).$$
In this paper, we consider this ordinary notion of subadditivity is of order $1$ and generalized the concept for any order $n$, where $n\in\N$. We establish  that $n^{th}$ square root of a $n^{th}$ order subadditive function possesses ordinary subadditivity. We also introduce the notion of approximately subadditive function and showed that it can be decomposed as the algebraic summation of a subadditive and a bounded function. Another important newly introduced concept is Periodical monotonicity. A function $f:I\to\R$ is said to be periodically monotone with a period $d>0$ if the following holds
\Eq{*}{
f(x)\leq f(y)\qquad\mbox{for all}\quad x,y\in I\qquad{with}\quad y-x\geq d.
}
One of the obtained results is that under a minimal assumption on $f$; this type of function can be decomposed as the sum of a monotone and a periodic function whose period is $d$. Towards the end of the paper, we discuss about star convexity. 

A function $f: I\to\R$ is said to be star-convex if there exists a point $p\in I$ such that for any $x\in I$ and for all $t\in [0,1]$; it satisfies either one of the following conditions.
$$
t(x,f(x)) +(1-t)(p,f(p))\in epi(f) \quad \mbox{or} \quad hypo(f).
$$
We studied the structural properties and showed relationship of it with star convex bodies. The motivation and other crucial details are discussed in the following section.
\end{abstract}
\maketitle
\section{Introduction}
This paper can be divided into three folds. At first  we consider $I$ to be a non negative interval and define notions of higher order and approximate subadditivity. Next, we consider the concept of distance and discuss about distance oriented monotonicity. At last, we generalized star-like convexity. We studied decomposition, sandwich type theorems and other functional properties of these classes of functions. At the beginning of each section a proper analysis of structural characteristics are also given.

Some crucial studies on subadditive functions were made in the  papers (cf.\ \cite{Mat93, Ros, Bin, Ham62, Meh, Lat}). Here studies related to measurable subadditive functions, Hermite-Hadamard type inequalities and some generalized forms for subadditive functions were conducted in depth.
In the direction of generalization of subadditivity our motivation comes from the theory of approximate monotonicity and convexity which has become an active field of research and thus many important contributions have been made. In our accepted paper of \cite{890}, a generalized form of subadditivity was discussed which was termed as $\Gamma$- property. A short description regarding background of $\Gamma$- property is given below.

The following two inequalities involved in each of the pairs are  equivalent to one another.  These assertions were established in the papers \cite{GosPal20} and \cite{890}  respectively.
\begin{enumerate}[(i)]
\item If $f:I\to\R$ is approximately monotone with the error function $\Phi:[0,\ell(I)]\to\R_+$ then
\Eq{*}{
f(x)\leq f(y)+\Phi(y-x)\quad & \mbox{for all}\quad x,y\in I\quad \mbox{with} \quad x\leq y\\
&\equiv \\
f(x)\leq f(y)+\Phi^{\sigma}(y-x)\quad & \mbox{for all}\quad x,y\in I \quad \mbox{with} \quad x\leq y.
}
\item While if $f:I\to\R$ is approximately convex with the error function $\Phi:[0,\ell(I)]\to\R_+$ then for any $x,y \in I$ and for all $t\in [0,1]$ the following will satisfy
\Eq{*}{
f(tx+(1-t)y)\leq tf(x)+(1-t)f(y)+& t\Phi((1-t)|y-x|)+(1-t)\Phi(t|y-x|)\\
&\equiv\\
f(tx+(1-t)y)\leq tf(x)+(1-t)f(y)+&t\Phi^{\gamma}((1-t)|y-x|)+(1-t)\Phi^{\gamma}(t|y-x|).
}
\end{enumerate}
It was proved that for the first equivalency result, $\Phi^{\sigma}$ is the largest subadditive minorant of $\Phi$. On the other hand, in the second assertion $\Phi^{\gamma}$ is the largest  minorant of $\Phi$ that possesses $\Gamma$- property. 

A function $\Phi: [0,\ell(I)]\to\R_+$ is said to be equipped with $\Gamma$-property if for any $x,y$ and $x+y\in [0,\ell(I)]$ provided $y>0$; it satisfies the following functional inequality 
\Eq{*}{
\Phi(x+y)\leq \Phi(x)+\dfrac{(x+y)^2-x^2}{y^2}\Phi(y)
.}
In the paper \cite{890}, several relationships between $\Phi^{\gamma}$ and $\Phi^{\sigma}$ were also studied. Such as $\sqrt{\Phi^{\gamma}}$ and $\dfrac{\Phi^{\gamma}(t)}{t}$ are subadditive functions in $[0,\ell(I)]$ and $]0,\ell(I)]$ respectively. One can easily observe that any ordinary subadditive function also satisfies $\Gamma$ property. In other words, we conclude that this function class contains the class of ordinary subaddative functions. Motivated from this and by definition, we rename $\Gamma$ property as second order subadditivity. We tried to generalize this newly introduced concept to any order $n\in\N$ and obtain some interesting results that can be later used in the theory of higher order convexity.

Let $I\subseteq \R_+$ be a non empty and non singleton interval. A function $\Phi: I\to \R_{+}$ is said to be subadditive of order $n\in \N$, if for any $x,y$ and $x+y\in I$ provided $y>0$, it satisfies the following functional inequality
\Eq{*}{
\Phi(x+y)\leq \Phi(x)+\frac{(x+y)^n-x^n}{y^n}\Phi(y).
} 

Based on the definition of $n$-order  subadditivity, we obtain some close relationships with ordinary subadditive function. We will see that, the 
$n^{th}$ square root of an $n$-order subadditive function is  simply subadditive. Besides, we also proved that the ratio of any $n$-order subadditive function to the $n^{th}$ power function is also going to be subadditive. Moreover, we also studied a functional inequality and establish the existence and uniqueness of its solution.
Next we define the notion of approximately subadditive functions.
 
A function $\Phi:I(\subseteq\R_+)\to\R_+$ is said to be approximately subadditive if there exists an $\varepsilon\geq0$ such that for any
arbitrary number of $u_1,\cdots u_n\in I$ assuring $u_1+\cdots u_n\in I$; satisfies following inequality
\Eq{*}{
\Phi(u_1+\cdots+u_n)\leq \Phi(u_1)+\cdots+\Phi(u_n)+\varepsilon.
}
We showed that a approximately subadditive function 
can be decomposed as the summation of a subadditive minorant  and a non negative function bounded above by $\varepsilon$.
In the next section we will go through  distance dependent monotonic functions. 

The notions of approximately monotone, sub-monotone functions were coined up in the papers (cf.\ \cite{Csa60,DanGeo04,ElsPeaRob96,MakPal12a,NgaPen07,Pal03a}) and have applications in nonsmooth and convex analysis and optimization theory, and also in the theory of functional equations and inequalities. In one of our previous papers \cite{GosPal20}, we introduced the concept of approximately monotone functions and broadly studied their structural, characteristic and functional properties. We now recall the terminology and notation of monotonicity introduced in that paper.

Let, $I$ be an non empty and non singleton open interval. A function $f:I\to\R$ is said to be \emph{$\Phi$-monotone} if, for any $x,y\in I$ with $x\leq y$ the following inequality holds 
\Eq{*}{
  f(x)\leq f(y)+\Phi(y-x),
}
where $\Phi: [0, \ell(I)]$ is termed as error function.

If this inequality is satisfied with the identically zero error function $\Phi$, then we can see that $f$ is \emph{monotone (increasing)}. Suppose, $\ell(I)\geq d$ and for any $x,y\in I$ with $y-x \geq d$; $\Phi(y-x)=0$ holds; then the function $f$ can be treated as distance dependent monotone(increasing) function.

Motivated by this, we introduce the following concept. Let $d>0$. A function $f: I\to \R$ is said to be increasing by a period $d$ (or $d$- periodically increasing) if the following  holds
$$f(x)\leq f(y) \qquad \mbox{ for all    }\quad x,y\in I\quad\mbox{with} \quad  y-x\geq d.$$
We start with some of the fundamental structural properties. For any given function $f$ we give a precise formula to obtain the largest $d$-periodically increasing minorant. We also show approximation of $d$-periodically increasing function by an ordinary monotone function.

Finally, we generalized the concept of star-like convexity introduced and characterized in the papers  \cite{Losonczi2020} and \cite{889} respectively. A function $f: I\to \R$ is said to be star-convex if there exists atleast an element $p\in I$ such that for any $x\in I$ and for all $t\in[0,1]$; the entire line segment joining $p$ and $x$, $t(x,f(x))+(1-t)(p,f(p))\in epi(f)$ or  $hypo(f)$. We showed that join of a convex/concave function with a another convex/concave; can also generate star-convex function in the elaborated domain . We ended our paper by studying the appearance of star convexity in appropriate combination epigraph and hypograph of star-convex function.

\section{On higher Order Subadditivity} 
Throughout this paper $\R_+$ will be the set of non negative numbers and for this section we will consider $I\subseteq \R_+$ be an non empty and non singular interval. A function $\Phi: I\to\R_{+}$ is said to be subadditive of order $n(\in\N)$ if, for any $x,y$ and $x+y\in I$ provided $y>0$, satisfies the following functional inequality
\Eq{B111}{
\Phi(x+y)\leq \Phi(x)+\frac{(x+y)^n-x^n}{y^n}\Phi(y).
} 
It is evident that if both $\Phi$ and $\Psi$ are two subadditive function of order $n;$  Then $\Phi+\Psi$ will also be $n$-order subadditive. Moreover, for any $a>0$; the function $(a\Phi)$ will also posses $n$-order subadditivity. Let $\F=\{\Phi_\gamma\mid\gamma\in \Gamma\}$ be a family of subadditive functions of order $n$ with a pointwise supremum $\Phi$. Then $\Phi$ is also going to be a $n$-order subadditive.  

\Thm{99999}{Let $\Phi: I\to\R_+$ be a subadditive function of order $n$. Then $\Phi$ is also be a subadditive function of order $n+1$.}

\begin{proof} 
For the proof, first we need the following inequality: For $u\geq0$,
\Eq{579}{
   (u+1)^n-u^n\leq (u+1)^{n+1}-u^{n+1}.
}
To establish this inequality, we proceed as follows
\Eq{*}{
(u+1)^{n+1}-u^{n+1}-(u+1)^n+u^n
&=u(u+1)^{n}-u^{n}(u-1)\\
&=u \bigg(\binom{n}{1} u^{n-1}+\cdots+1\bigg)+u^{n}\geq 0
.}
$\Phi$ is $n$- order subadditive; it will satisfy \eq{B111}.
Together with it, by replacing $u$ as $\dfrac{x}{y}$ in \eq{579} ; we obtain
 
\Eq{*}{
\Phi(x+y)\leq \Phi(x)+\frac{(x+y)^n-x^n}{y^n}\Phi(y)
\leq \Phi(x)+\frac{(x+y)^{n+1}-x^{n+1}}{y^{n+1}}\Phi(y).
}
This yields that $\Phi$ is also of $(n+1)$-order subadditive and completes the proof.
\end{proof}
The following corollary establishes the linkage of power function with the $n$-order subadditivity. It is also evident that even super additive functions can also possess higher order subadditivity.

\Cor{1113}{The power function $x\to x^{p}$ $(x\geq 0)$ is subadditive of order $n$, where $n=\lceil p\rceil$.}
\begin{proof}
Substituting $\Phi(x)=x^{n}$ in $\eq{B111}$, we obtain
\Eq{776}{
(x+y)^{n}= x^{n}+\dfrac{(x+y)^{n}-x^{n}}{y^{n}}y^{n}.
}
Which shows that the power function $t^{p}$ will be subadditive of order $n$ if $p=n$.
By \thm{99999}, the statement is also valid for $p\leq n-1$.

To prove the corollary, we need to consider the only remaining case; i.e when $p\in ]n-1,n[$. Let, $p\in ]n-1,n[$ is  arbitrary. Then there exists an element $q\in]0,1[$ such that $p+q=n$. From \eq{776} we obtain
\Eq{*}{
(x+y)^{p}&=\dfrac{1}{(x+y)^q} \bigg[x^{n}+\dfrac{(x+y)^{n}-x^{n}}{y^{n}}y^{n}\bigg]\\
&\leq \dfrac{x^{n}}{x^q}+\bigg[\dfrac{(x+y)^{n}-x^{n}}{y^{n}}\bigg]\dfrac{y^n}{y^q}\\
&= x^p+\bigg[\dfrac{(x+y)^{n}-x^{n}}{y^{n}}\bigg]y^p.}
This completes the proof.
\end{proof}
The following two theorems yields the close linkage of ordinary subadditivity with the higher order subadditive function. 
\Thm{B2}{Let $\Phi: I\to \R_{+}$ is $n$-subadditive. Then $\sqrt[n]{\Phi}$ is subadditive.}
\begin{proof}
For the proof, first we need to establish the below mentioned inequality.

For $0\leq u\leq v$, the following holds
\Eq{5799}{
   (u+1)^n-u^n\leq (v+1)^{n}-v^{n}.
}
To establish this inequality, we simply observe that
\Eq{*}{
(v+1)^{n}-v^{n}-(u+1)^n+u^n
&=\sum_{i=1}^{n-1}{n \choose i}\Big(v^{n-i}-u^{n-i}\Big)\geq 0 .}
To prove the assertion, Let $x,y\in I$ such that $y>0$ and $x+y\in I$. without loss of generality, we assume that 
\Eq{0023}{\dfrac{x}{y}\leq \dfrac{\sqrt[n]{\Phi(x)}}{\sqrt[n]{\Phi(y)}}.}
Since $\Phi$ is subadditive of order $n$; therefore it will satisfy $\eq{B111}$. By using \eq{0023} and utilizing \eq{5799} there, we obtain
\Eq{*}{
\Phi(x+y)&\leq \Phi(x)+\frac{(x+y)^n-x^n}{y^n}\Phi(y)
\\
&\leq \Phi(x)+\Bigg(\frac{\sqrt[n]{\Phi(x)}}{\sqrt[n]{\Phi(y)}}
+1\Bigg)^{n}\Phi(y)-\Bigg(\dfrac{\sqrt[n]{\Phi(x)}}{\sqrt[n]{\Phi(y)}}\Bigg)^{n}\Phi(y)\\
&=\Big(\sqrt[n]{\Phi(x)}+\sqrt[n]{\Phi(y)}\Big)^{n}.}
Upon taking the $n^{th}$ square root on the both side of the above inequality we obtain
\Eq{209}{
\sqrt[n]{\Phi(x+y)}\leq \sqrt[n]{\Phi(x)}+\sqrt[n]{\Phi(y)};}

Similarly in the reverse case of \eq{0023}, we can proceed in the similar way by interchanging the roles of $x,y$ in \eq{B111}. Finally we will obtain the same inequality as \eq{209}; which yields that $\sqrt[n]{\Phi}$ is subadditive. It completes our proof.
\end{proof}
\Thm{T3}{
$\Phi: I\to\R_+$ be a subadditive function of order $n$. Then $\dfrac{\Phi(t)}{t^{n}}$ will be subadditive in $I\cap]0,\infty[$.
}
\begin{proof}
Since $\Phi$ be a subadditive function of order $n$, then  it will satisfy \eq{B111}. 
By interchanging the roles of $x$ and $y$ we get
\Eq{B211}{
\Phi(x+y)\leq \Phi(y)+\dfrac{(x+y)^{n}-y^{n}}{x^{n}}\Phi(x).}
Summing up \eq{B111} and \eq{B211} side by side we arrive at
\Eq{*}{
2\Phi(x+y)&\leq \big(x^n-y^n+(x+y)^n\big)\dfrac{\Phi(x)}{x^n}+\big(y^n-x^n+(x+y)^n\big)\dfrac{\Phi(x)}{x^n}\\
&\leq 2(x+y)^n\bigg[\dfrac{\Phi(x)}{x^n}+\dfrac{\Phi(y)}{y^n}\bigg]
.}
Dividing both side of the above inequality by $2(x+y)^n$, we obtain
\Eq{*}{
\dfrac{\Phi(x+y)}{(x+y)^n}\leq \dfrac{\Phi(x)}{x^n}+\dfrac{\Phi(y)}{y^n}.
}
This shows the subadditivity of $\dfrac{\Phi(t)}{t^{n}}$ in $I\cap]0,\infty[$ and hence completes the proof.
\end{proof}
The next proposition establishes a weaker form $n^{th}$ order subadditivity by omitting the binomial combination of $x$ and $y$ to the right hand side of \eq{B111}.
\Prp{6777} {Let $\Phi: I\to\R$ be a $n$-order subadditive function. Then the following inequality holds
\Eq{998}{
\Phi(x+y)\leq \max\{\Phi(x)+(2^n-1)\Phi(y),(2^n-1)\Phi(x)+\Phi(y)\}
}}
\begin{proof}
To prove the statement without loss of generality assume $x\leq y$. Then, by definition of $n$-order subadditivity of $\Phi$, we have
\Eq{*}{
\Phi(x+y)&\leq \Phi(x)+\frac{(x+y)^n-x^n}{y^n}\Phi(y)\\
&= \Phi(x)+\bigg[\binom{n}{1}\bigg(\frac{x}{y}\bigg)^{n-1}+\cdots+\binom{n}{n-1}\bigg(\frac{x}{y}\bigg)+1\bigg]\Phi(y)\\
&\leq \Phi(x)+\bigg[\binom{n}{1}+\cdots+\binom{n}{n-1}+1\bigg]\Phi(y)\\
&=\Phi(x)+(2^n-1)\Phi(y).}
Conversely if $y\leq x$ hold then by using \eq{B211}, we can show
\Eq{*}{
\Phi(x+y)\leq (2^n-1)\Phi(x)+\Phi(y).}
Combining these two above inequalities we obtain \eq{998} and it proves the statement.
\end{proof}
Now we will discuss about a functional equation and we will show that it has one and only one solution form.
\Prp{01}{Let $I\subseteq \R_+\setminus\{0\}$ be a non empty and non singleton interval and $n(\in \N$) be any number such that $n\geq2.$ The function  $\Phi: I\to \R$ satisfies the following functional equality 
\Eq{B11111}{
\Phi(x)+\frac{(x+y)^n-x^n}{y^n}\Phi(y)=\Phi(y)+\frac{(x+y)^n-y^n}{x^n}\Phi(x).
}
Then the only possible solution  will be in the form $$\Phi(x)=cx^n\qquad\mbox{where} \quad c\in \R.$$}
\begin{proof}
Rearranging the terms of \eq{B11111} we obtain
\Eq{*}{
\bigg[(x+y)^{n}-\Big(x^{n}+y^{n}\Big)\bigg]\bigg[\dfrac{\Phi(x)}{x^{n}}-\dfrac{\Phi(y)}{y^{n}}\bigg]=0.}
 
It is evident that for $n\geq 2$, the expression $(x+y)^{n}-\big(x^{n}+y^{n}\big)>0.$
Therefore to satisfy the above equality the following must hold
\Eq{*}{\dfrac{\Phi(x)}{x^n}=\dfrac{\Phi(y)}{y^n} \qquad (\mbox{for all } x,y\in I\,\, \mbox{ and for any } n\in \N).}
Which shows that $\dfrac{\Phi(x)}{x^n}$ is a constant function on $I$ and hence the solution of the functional inequality \eq{B11111} will be in the form $\Phi(x)=cx^n$, where $c\in \R$.

\end{proof}

A function $\Phi: I\to\R_+$ is said to be approximately subadditive if there exists an $\varepsilon>0$ such that for any number of $u_1,\cdots u_n\in I$ provided $u_1+\cdots u_n\in I$; satisfies the following inequality
\Eq{*}{
\Phi(u_1+\cdots+u_n)\leq \Phi(u_1)+\cdots+\Phi(u_n)+\epsilon.
}
The following theorem gives a possible decomposition of approximately subadditive function. The theorem extensively  use Proposition 3.3 of \cite{GosPal20}.
\Thm{sa}{$\Phi: I\to\R_+$ is approximately 
subadditive if and only if it can be decomposed as the sum of $\Phi^{\sigma}$ and $h$; where $\Phi^{\sigma}: I\to\R_+$ is the largest subadditive minorant of $\Phi$ and h is a non negative function satisfying $h\leq\epsilon$. Additionally if $\Phi$ is increasing then $h$ is a non negative function of bounded variation.}
\begin{proof}
Since $0\leq\Phi$; it ensures that there exists atleast one subadditive minorant for $\Phi$.

First, we assume that $\Phi^{\sigma}: I\to\R_+$ is the greatest subadditive minorant of $\Phi$ and $h\leq\epsilon$. Then for any number of $u_1,\cdots u_n\in I$ satisfying $u_1+\cdots +u_n\in I$, the following two inequalities holds
\Eq{*}{
\Phi^{\sigma}(u_1+\cdots+u_n)&\leq \Phi^{\sigma}(u_1)+\cdots+\Phi^{\sigma}(u_n)\\
&\mbox{and}\\
h(u_1+&\cdots+u_n)\leq \epsilon.
}
Summing up these two inequalities side by side we obtain
\Eq{*}{
\Phi(u_1+\cdots+u_n)&=\Phi^{\sigma}(u_1+\cdots+u_n)+h(u_1+\cdots+u_n)\\
&\leq \Phi^{\sigma}(u_1)+\cdots+\Phi^{\sigma}(u_n)+\epsilon\\
&\leq \Phi(u_1)+\cdots+\Phi(u_n)+\epsilon.
}
This establishes that $\Phi$ is approximately subadditive.
Conversely, suppose $\Phi$ is approximately subadditive.
Now, we define the function $\Phi^\sigma:I\to\R_+$ by
\Eq{*}{
\Phi^\sigma(u):=\inf\big\{\Phi(u_1)+\dots+\Phi(u_n)\mid
  n\in\N, u_1,\dots,u_n\in I\colon u_1+\dots+u_n=u\big\}.
}
We will show that $\Phi^\sigma$ is the largest subadditive function which satisfies the inequality $\Phi^\sigma\leq\Phi$ on $I$.
First we are going to prove the subadditivity of $\Phi^\sigma$. Let $u,v$ and $u+v\in I$. Let $\delta>0$ be arbitrary. Then there exist $n,m \in\N$ and $u_1,\dots,u_n,v_1,\dots,v_m\in I$ such that 
\Eq{v}
{u=\sum_{i=1}^{n} u_i,\qquad v=\sum_{j=1}^{m} v_j, \qquad
\sum_{i=1}^{n}\Phi(u_i)<\Phi^\sigma(u)+\frac{\delta}{2} 
\qquad\text{and}\qquad
\sum_{j=1}^{m}\Phi(v_j)<\Phi^\sigma(v)+\frac{\delta}{2}.
}
We have that $u+v=\sum_{i=1}^{n} u_i+\sum_{j=1}^{m} v_j$. Therefore, by the definition of $\Phi^\sigma$ and by the last two inequalities in \eq{v}, we get
\Eq{*}{
  \Phi^\sigma(u+v)\leq \sum_{i=1}^{n}\Phi(u_i) + \sum_{j=1}^{m}\Phi(v_j)
  <\Phi^\sigma(u)+\Phi^\sigma(v)+\delta.
}
Since $\delta$ is an arbitrary positive number, we conclude that $\Phi^\sigma(u+v)\leq \Phi^\sigma(u)+\Phi^\sigma(v)$, which shows the subadditivity of $\Phi^\sigma$.

 By taking $n=1$, $u_1=u$ in the definition of $\Phi^\sigma(u)$, we can see that $\Phi^\sigma(u)\leq\Phi(u)$.

 Now assume that $\Psi:I \to\R_+$ is a subadditive function such that $\Psi\leq\Phi$ holds on $I$. To show that $\Psi\leq\Phi^\sigma$, let $u\in I$ and $\delta>0$ be arbitrary. Then 
there exist $n \in\N$ and $u_1,\dots,u_n\in I$ such that conditions for $u$ in \eq{v} are hold.
Then, due to the subadditivity of $\Psi$,
\Eq{*}{
  \Psi(u)\leq \Psi(u_1)+\dots+\Psi(u_n)
  \leq \Phi(u_1)+\dots+\Phi(u_n)<\Phi^\sigma(u)+\dfrac{\delta}{2}.
}
By the arbitrariness of $\delta>0$, the inequality $\Psi(u)\leq\Phi^\sigma(u)$ follows for all $u\in I$. It establishes that $\Phi^{\sigma}$ is the largest subadditive function satisfying $\Phi^{\sigma}\leq\Phi.$

Now it will be enough to show $h:=\Phi-\Phi^{\sigma}\leq \epsilon$. Let $u\in I$ be arbitrary. Then there exits a non negative partition of $u$; $ \{u_1,\cdots u_m\}$ $(u_i's\in I)$ satisfying $\sum_{i=1}^{m}\Phi(u_i)<\Phi^{\sigma}(u)+\delta$, where $\delta$ is a negligibly small positive number. Now using the definition of $\Phi^{\sigma}$, we have the following oservation
\Eq{*}{
h(u)&=\Phi(u)-\Phi^{\sigma}(u)\\
 &\leq  \Big({\Phi(u_1)+\dots+\Phi(u_m)+\epsilon\Big)-\Phi^{\sigma}(u)}<\varepsilon+\delta.}
Taking $\delta\to 0$, we finally obtain $0\leq h\leq\varepsilon$ and it completes the decomposition part.

To verify the last assertion, let $u,v\in I$ with $v<u$. Let $\delta>0$ be arbitrary. Then there exist $n \in\N$ and $u_1,\dots,u_n\in I$ such that conditions for $u$ in \eq{v} is satisfied. Define $v_i:=\frac{v}{u}u_i$. Then $v_i\leq u_i$ and hence $\Phi(v_i)\leq\Phi(u_i)$. On the other hand, $v_1+\dots+v_n=\frac{v}{u}(u_1+\dots+u_n)=v$, which implies that
\Eq{*}{
  \Phi^\sigma(v)
  \leq \Phi(v_1)+\dots+\Phi(v_n) 
  \leq \Phi(u_1)+\dots+\Phi(u_n)
  <\Phi^\sigma(u)+\delta.
}
Passing the limit $\delta\to0$, we arrive at the inequality $\Phi^\sigma(v)\leq\Phi^\sigma(u)$, which proves the increasingness of $\Phi^\sigma$. Now by Jordan Decomposition Theorem  $h$ can be expressed as a difference of two monotonic function $\Phi$ and $\Phi^\sigma$. And hence $h$ is a function of bounded variation.
This completes the proof.
\end{proof}
\section{On Periodically increasing functions}
Let $d>0$ and $I\subseteq\R$ be a non empty, non singleton interval with $\ell (I)\geq d$. A function $f: I\to \R$ is said to be increasing by a period $d$ if for any $x,y\in I$, with $y-x\geq d$, $f(x)\leq f(y)$ holds.
For the rest of the section, we will call such function as $d$-periodically increasing function. 

One can easily observe that the class of $d$ periodically increasing function is a convex cone as it is closed under addition and multiplication with non negative scalars. Moreover if $f$ is $d$-periodically increasing and non negative then  $f^{n}$ is also periodically increasing with the same period $d$ for any $n(\in \N)$. The class of $d$-periodically increasing functions  is also closed under pointwise limit operations.

A function $f:I\to\R$ is said to have extended height $\H_{x_{d}}(f)$ in the interval $I_{d}=[x,x+d]\cap I\subseteq I$, which is expressible the following way
\Eq{*}{\H_{x_{d}}(f):=\sup_{u,v\in I_d}|f(u)-f(v)|.}
The supremum of all such $\H_{x_{d}}(f)$ will be denoted by $\H_d(f)$ and can be described as 
\Eq{*}{\H_d(f)=\sup_{x\in I}\H_{x_{d}}(f)=\sup_{I_{d}\subseteq I}\Bigg(\sup_{u,v\in I_d}|f(u)-f(v)|\Bigg).}
Finally for a function $f:I\to\R$; the height $\H(f)$ defined as follows
\Eq{*}{\H(f):=\sup_{u,v\in I}|f(u)-f(v)|.}
The following theorem provides a decomposition of a $d$-periodically increasing function under minimal assumptions.

\Thm{C1}{ Let $g: I\to\R$ be an increasing and  $k:I\to \R$ is a bounded function such that for any $x\in I$, $\H_{x_{d}}(g)\geq \H(k)$ holds. Then both $g+k$ and $g-k$ will be $d$-periodically increasing. Conversely let $\ell(I)>2d$. Then if $f$ is a $d$-periodically increasing function satisfying  $f(x+nd)-f(x+(n-1)d)=l>0,\, (n\in\N)$ for all $x\in I$ with $[x+nd, x+(n-1)d]\cap I\subseteq I$; then $f$  can be written as the sum of an increasing and a periodic function.}
\begin{proof}
To prove the first part of the theorem, we assume $x,y\in I$ with $x<x+d\leq y$.
By the definition of $g$ and $h$ and using the mentioned conditions, we get the following inequalities 
\Eq{*}{
k(x)-k(y)\leq \H(k)\leq\H_{x_{d}}(g)&=g(x+d)-g(x)\leq g(y)-g(x)\\
&\mbox{and}\\
k(y)-k(x)\leq \H(k)\leq\H_{x_{d}}(g)&=g(x+d)-g(x)\leq g(y)-g(x)\\
}
Rearranging these two inequalities we obtain
\Eq{*}{(g+k)(x)\leq (g+k)(y)\quad\mbox{and}\quad(g-k)(x)\leq (g-k)(y).
}
Which yields that both $g+k$ and $g-k$ are $d$ periodically increasing functions.

Now to prove the converse assertion, let $f$ be a $d$-periodically increasing function such that all the assumed(pre-mentioned) conditions holds. 
Now we define $g: I\to \R$ as 
$$g(x):=\inf_{x\leq t}f(t).$$
Clearly $g\leq f$.
Let $x,y\in I$ with $x\leq y$. Then by the following inequality we see that $g$ is a monotone(increasing) function in $I$
\Eq{334}{
g(x)=\inf_{x\leq t}f(t)\leq \inf_{y\leq t}f(t)\leq g(y).}
  
Let $x\in I$ be an arbitrary point. Since, $\ell(I)>2d$; therefore atleast either $x+d$ or $x-d\in I$. If not then
\Eq{*}{
x-d\leq \inf(I)\qquad{and} \qquad x+d\geq \sup(I)
.}
This implies $\ell(I)=\sup(I)-\inf(I)\leq(x+d)-(x-d)\leq 2d$.
Which results in a contradiction.

Now without loss of generality we assume that $x+d\in I$ 
and consider $h:=f-g$. It will be sufficient to show $h$ is $d$-periodic.  we have
\Eq{789}{
h(x+d)-h(x)=&\Big(f(x+d)- \inf_{x+d\leq t}f(t)\Big)-\Big(f(x)-\inf_{x\leq t}f(t)\Big)\\
=& \Big(f(x+d)-f(x)\Big)-\Big(\inf_{t\in [x+d,x+2d]\cap I}f(t)-\inf_{t\in[x,x+d]}f(t)\Big)\\
=& l-\Big(\inf_{t\in [x+d,x+2d]\cap I}f(t)-\inf_{t\in[x,x+d]}f(t)\Big)
.}
We claim that, if 
 $\inf_{t\in [x+d,x+2d]\cap I}f(t)=f(p)$ ($p\in [x+d,x+2d]$) this implies $\inf_{t\in [x,x+d]}f(t)=f(p-d)$. And hence by our assumed condition on $f$ the following inequality holds
$$\inf_{t\in[x+d,x+2d]\cap I} f(t)-\inf_{t\in[x,x+d]}f(t)=l.$$
If possible, let  $\inf_{t\in [x,x+d]}f(t)=f(q)$ such that $q\neq p-d$.
 Then the following inequality holds
\Eq{*}{
f(q)< f(p-d)<f(p)\leq f(q+d)
.}
This inequality yields $f(p)-f(p-d)<f(q+d)-f(q).$ It contradicts to our original assumed condition on $f$ that, for any $x\in I$ with $x+nd,x+(n-1)d\in I\,(n\in\N)$, the difference $f(x+nd)-f(x+(n-1)d)=l$.

Therefore $\inf_{t\in [x,x+d]}f(t)=f(p-d)$ and hence from \eq{789} we obtain
\Eq{*}{h(x+d)-h(x)=l-l=0.}
This shows $h$ is $d$-periodic and completes the proof.
\end{proof}
\Thm{C2}{Let, $f: I\to \R$ is bounded from below, then
the function $\widetilde{f}: I\to \R$ defined by
\Eq{*}{
\widetilde{f}(x):=\min\Big\{f(x), \inf_{ [x+d,\infty[\cap I}f(t)\Big\}\qquad \mbox{for all} \quad (x\in I)}
 will be the greatest $d$- periodically increasing minorant of $f$.}
\begin{proof}
To prove this theorem, let $x$ and $y\in I$ are two arbitrary elements with $y-x\geq d$. Now if $\widetilde{f}(x)=f(x)$, we have
\Eq{*}{
\widetilde{f}(x)=f(x)\leq \inf_{[x+d,\infty[\cap I}f(t)&= \min\Big\{f(y), \inf_{[x+d,\infty[\cap I}f(t)\Big\}\\
&\leq\min\Big\{f(y), \inf_{[y+d,\infty[\cap I}f(t)\Big\}=\widetilde{f}(y)
.}
On the other hand, if $\widetilde{f}(x)= \inf_{[x+d,\infty[\cap I}$, we obtain 
\Eq{*}{
\widetilde{f}(x)= \inf_{[x+d,\infty[\cap I}f(t)\leq \min\Big\{f(y), \inf_{[y+d,\infty[\cap I}f(t)\Big\}=\widetilde{f}(y)
.}
This yields that $\widetilde{f}$ is $d$-periodically increasing.

To show the last part of the theorem, we assume that  there exists a $d$-periodic increasing function $g:I\to\R$ such that $\widetilde{f}\leq g\leq f$ holds. Therefore, for any $x\in I$; and for all $t\geq x+d$; will satisfy the following function inequalities  
\Eq{*}{
g(x)&\leq f(x)\\ 
&\mbox{and}\\
g(x)\leq g(t)\leq f(t) \quad \mbox{  which  } & \mbox{  implies  }\quad g(x)\leq \inf_{x+d\leq t}f(t).
}
Combining the above two inequalities we have
\Eq{*}
{g(x)\leq \min\Big\{f(x),\inf_{[x+d,\infty[\cap I}f(t)\Big\}=\widetilde{f}(x).}
This establishes that $\widetilde{f}$ is the largest $d$-periodically increasing minorant of $f$.
\end{proof}
\Thm{212}{
Let $f$ be a $d$ periodically increasing function. Then there exists an increasing function $\widehat{f}$ such that $||f-\widehat{f}||\leq \dfrac{\H_{d}(I)}{2}$
}
\begin{proof}
Let $f$  be a $d$-periodically increasing function. Then we can define the largest and smallest monotone (increasing) minorant and majorant $f_{*}, f^{*}: I\to \R$ of $f$ as follows:
\Eq{*}{
f_{*}(x):=\inf_{x\leq t}f(t)\qquad \mbox{and}\qquad f^{*}(x):=\sup_{t\leq x}f(t).
}
First we will show that  $f_{*}$ is the largest monotone(increasing) function satisfying $f_*\leq f$. We assume $x,y\in I$ with $x<y$. Then substituting $f_*$, instead of $g$ in \eq{334} we obtain the monotonicity. If possible, let $g: I\to\R$ be the largest monotone function satisfying $f_*<g\leq f.$
Then we have the following inequality
\Eq{*}{
g(x)\leq g(t) \leq f(t)\qquad \mbox{for all}\quad x\leq t.
}
Taking infimum to the right hand side of the above inequality we obtain 
\Eq{*}{
g(x)\leq \inf_{x\leq t}f(t)=f_*(x)\qquad \mbox{for all}\quad x\in I.
}
which contradicts our assumption that $f_*<g$. Therefore, $f_{*}$ is the largest monotone (increasing) minorant of $f.$ 

Similarly, we can also show that $f^*$ is the smallest monotone function satisfying of $f\leq f^*$.

Now we can define the increasing function 
 $\widehat{f}=:\dfrac{f^{*}+f_{*}}{2}$ and calculate upper bounds for $f-\widehat{f}$ and $\widehat{f}-f$ respectively as follows. Let, $x\in I$ be arbitrary
\Eq{*}{
(f-\widehat{f})(x)&=\bigg(\dfrac{f-f_{*}}{2}\bigg)(x)+\bigg(\dfrac{f-f^{*}}{2}\bigg)(x)\\
&\leq \dfrac{1}{2}\Big(f(x)-f_{*}(x)\Big)\\
&=\dfrac{1}{2}\Big(f(x)-\inf_{x\leq t}f(t)\Big)\\
&=\dfrac{1}{2}\Big(f(x)-\inf_{[x,x+d]\cap I}f(t)\Big)\\
&\leq \dfrac{\H_{x_{d}}(I)}{2}\leq \dfrac{\H_d(I)}{2}
.}
Therefore $f-\widehat{f}\leq \dfrac{\H_d(I)}{2}$.
Similarly, we can also show that $\widehat{f}-f\leq \dfrac{\H_d(I)}{2}$. These two inequalities together prove the statement.
\end{proof}
\section{On Star-Convex Functions}
Before proceeding, at first we need to recall the following three terminologies related to any function $f:I\to\R$; namely graph, epigraph and hypograph of the function. We will extensively use these concepts in the definition and to formulate various geometric characterization of star-convex function.
\Eq{*}{
gr(f):=\big\{(x,t)\in I\times\R\mid f(x)=t\big\}
.}
\Eq{*}{
epi(f):=\big\{(x,t)\in I\times\R\mid f(x)\leq t\big\}
.}
\Eq{*}{
hypo(f):=\big\{(x,t)\in I\times\R\mid f(x)\geq t\big\}
.}
A function $f: I\to \R$ is said to be star-convex if there exists atleast an element $p\in I$ such that for any arbitrary $x\in I$ and for all $t\in[0,1]$ either one of the following inequality holds
\Eq{*}{
f(tx+(1-t)p)\leq tf(x)+(1-t)f(p)\qquad \mbox{or}\qquad  f(tx+(1-t)p)\geq tf(x)+(1-t)f(p).
}
In other words for any arbitrary $x\in I$, the line segment joining the two points $x$ and $p$; will entirely lies either in $epi(f)$ or $hypo(f)$. That is for all $t\in[0,1]$ either $(tx+(1-t)p, tf(x)+(1-t)f(p))\in epi(f)$ or $hypo(f)$. The point $(p,f(p))\in gr(f)$ will be called the center of the function $f$. And the set of all such points $p$ is termed as central set. We will denote the central set of $gr(f)$ by $\C_{f}(I)$. Based on this definition it is evident that if $f$ is convex(concave) then $f$ will naturally be a star-convex with  $\C_{f}(I)=gr(f)$. Also the class of star-convex functions will be a centrally symmetric convex cone as it is closed under algebraic summation,  multiplication with non negative scalars and both $f$ and $-f$ are belong to the same class of function. 

The following proposition describes the structural composition of some specific types of star-convex functions. 
\Prp{223}{Let $f: I\to \R$ be a continuous function such that one of the following condition holds.
\begin{enumerate}[(i)]
\item  $f\Big \arrowvert_{]-\infty, p]}$ and  $f\Big \arrowvert_{[p,\infty]\cap I}$ are both convex.
\item  $f\Big \arrowvert_{]-\infty, p]}$  and  $f\Big \arrowvert_{[p,\infty]\cap I}$ are both concave.
\item $f\Big \arrowvert_{]-\infty, p]}$  is convex and   $f\Big \arrowvert_{[p,\infty]\cap I}$ is concave.
\item $f\Big \arrowvert_{]-\infty, p]}$  is concave and  $f\Big \arrowvert_{[p,\infty]\cap I}$ is convex.
\end{enumerate}
}
\begin{proof}
We will establish the proposition only for the condition $(i)$ as for all the other remaining cases, the proofs will be analogous.

Let, $f$ satisfies condition $(i)$. Therefore $(p,f(p)$ will be the common center for both $f\Big \arrowvert_{]-\infty, p]}$  and  $f\Big \arrowvert_{[p,\infty]\cap I}$. Which implies $(p,f(p)\in \C_f(I)$. This proves that $f$ is star convex.
\end{proof} 
A set $ X\subseteq \R^{2}$  is said to be a star-convex , if there exists atleast an $x_{0}\in X$ such that for all $x\in X$ the line segment from $x_{0}$ to $x$ lies entirely in $X$. The collection of all such $x_{0}\in X$ is called the central set of $X$. In other words, for a star-convex set $X$, the central set $\C$ is defined as follows:
$$\C=\{x_{0}\in X: \mbox{ for all } x\in X\mbox{ and } t\in[0,1]\mbox{, }tx+(1-t)x_{0}\in X\}.$$
If $X$ is convex or concave then $\C=X$.
The following proposition link the star-convex set with the star-convex function.
\Thm{8.112}{Let $f: I\to\R$ be a continuous star-convex function with $(p,f(p))\in C_{f}(I)$. Then the following assertions can be made.
\begin{enumerate}[(i)]
\item if condition $(i)$ of \prp{223} holds, then epi(f) is star-convex set.
\item if condition $(ii)$ of \prp{223} holds, then hypo(f) is star-convex set.
\item if condition $(iii)$ of \prp{223} holds, then $ \bigg(epi(f)\Big\arrowvert_{]-\infty,p]\cap I}\bigg)\cup \bigg( hypo(f)\Big\arrowvert_{ [p,\infty[\cap I} \bigg)$ is a star-convex set.
\item if condition $(iv)$ of \prp{223} holds,, then $ \bigg(hypo(f)\Big\arrowvert_{]-\infty,p]\cap I}\bigg)\cup \bigg( epi(f)\Big\arrowvert_{ [p,\infty[\cap I} \bigg)$ is a star-convex set.
\end{enumerate} }
\begin{proof}
Let $f: I\to\R$ is a star-convex function.
We will only show the validity for assertion $(i)$ as the remaining assertions will follow analogously. 

Assuming condition $(i)$ of \prp{223} on $f$ , we have that both of the sets 

$epi(f)\Big\arrowvert_{]-\infty,p]\cap I}$ and $epi(f)\Big\arrowvert_{[p,\infty[\cap I}$ are convex. Which implies $(p,f(p)$ will be a central point for $\Bigg(epi(f)\Big\arrowvert_{]-\infty,p]\cap I}\Bigg)\cup \Bigg(epi(f)\Big\arrowvert_{[p,\infty[\cap I}\Bigg)$.
This implies $(p,f(p))\in \C(epi(f)).$ This shows that $epi(f)$ is a star-convex set and completes our proof.
\end{proof}

\end{document}